\documentclass{amsart}
\usepackage{amsmath,amssymb,amscd}

\newcommand{\NI}{{\noindent}}
\newcommand{\QED}{\hfill$\Box$\medskip}
\newcommand{\Ham}{{\rm Ham}}
\newcommand{\Symp}{{\rm Symp}}
\newcommand{\R}{\mathbb{R}}
\newcommand{\Z}{\mathbb{Z}}
\newcommand{\T}{\mathbb{T}}
\newcommand{\A}{\mathcal{A}}
\newcommand{\M}{\mathcal{M}}

\numberwithin{equation}{section}

\newtheorem{theorem}{Theorem}[section]
\newtheorem{cor}[theorem]{Corollary}

\newtheorem{definition}[theorem]{Definition}

\newtheorem{remark}[theorem]{Remark}
\newtheorem{rmk}[theorem]{Remark}
\newtheorem{lemma}[theorem]{Lemma}

\newtheorem{guess}[theorem]{Conjecture}

\begin{document}

\title{A nonextension result on the spectral metric}
\author{Zhigang Han}
\address{Department of Mathematics and Statistics, University of Massachusetts, Amherst, MA 01003-9305, USA}
 \email{han@math.umass.edu}
\date{May 06, 2008}

\begin{abstract}
The spectral metric, defined by Schwarz and Oh using
Floer-theoretical method, is a bi-invariant metric on the
Hamiltonian diffeomorphism group. We show in this note that for
certain symplectic manifolds, this metric can not be extended to a
bi-invariant metric on the full group of symplectomorphisms. We also
study the bounded isometry conjecture of Lalonde and Polterovich in
the context of the spectral metric. In particular, we show that the
conjecture holds for the torus with all linear symplectic forms.
\end{abstract}

\maketitle

\section{Introduction}\label{introduction}
For a closed symplectic manifold $(M, \omega)$, the group $\Symp(M,
\omega)$ of all symplectomorphisms has two important subgroups,
namely, its identity component $\Symp_0(M, \omega)$ and the subgroup
$\Ham(M, \omega)$ of all Hamiltonian diffeomorphisms. These two
subgroups coincide when the first cohomology group $H^1(M, \R)$ of
$M$ vanishes, and they behave quite differently when $H^1(M, \R)$
does not vanish. For instance, a Hamiltonian diffeomorphism always
has a fixed point according to the Arnold conjecture, but a
symplectomorphism does not necessarily have one.

Another example is that the group $\Ham(M, \omega)$ admits several
interesting bi-invariant metrics, namely, the famous Hofer metric
discovered by Hofer \cite{Ho} and the spectral metric defined by
Schwarz \cite{Sch} and Oh \cite{Oh2} using Floer-theoretical method,
which we will describe in Section \ref{tsn}. For the group
$\Symp_0(M, \omega)$, although it also admits bi-invariant metrics
(cf. Prop 1.2.A \cite{LP} and Prop 3.10 \cite{Han2}), these metrics
are not quite naturally defined.

One way to construct interesting metrics on $\Symp_0(M, \omega)$ is
to extend the existing metrics such as the Hofer metric and the
spectral metric on $\Ham(M, \omega)$. In fact, such extensions are
constructed in \cite{Ban2, BD} for certain symplectic manifolds.
While the metrics on $\Ham(M, \omega)$ are bi-invariant, they only
extend to right-invariant metrics on $\Symp_0(M, \omega)$. In the
opposite direction, it is proven in \cite{Han2} that the Hofer
metric does not extend bi-invariantly to $\Symp_0(M, \omega)$ for
certain symplectic manifolds. In fact, it is conjectured that this
is true in general.

In this note, we focus on the spectral metric. In particular, using
ideas of Lalonde-Polterovich \cite{LP} and Schwarz \cite{Sch}, we
prove the following result. A submanifold $L \hookrightarrow M$ is
incompressible if the inclusion induces an injection on the
fundamental groups.

\begin{theorem}\label{thm1}
Let $(M, \omega)$ be closed and symplectically aspherical manifold,
and $L \hookrightarrow (M, \omega)$ be an incompressible Lagrangian
torus. If there exists some $\phi \in \Symp_0(M,\omega)$ such that
$\phi(L) \cap L =\emptyset$, then the spectral metric does not
extend to a bi-invariant metric on $\Symp_0(M,\omega)$.
\end{theorem}

\begin{rmk} \rm
The assumption that $M$ is symplectically aspherical can be removed,
as the spectral metric is defined for general closed symplectic
manifolds. However, we decide to keep it here since we only describe
the definition of the spectral metric for symplectically aspherical
manifolds in this paper.
\end{rmk}

The following corollary is immediate. Here $(\T^{2n},\omega_0)$
stands for the torus with the standard form.

\begin{cor}
The spectral metric on $\Ham(\T^{2n},\omega_0)$ does not extend to a
bi-invariant metric on $\Symp_0(\T^{2n},\omega_0)$.
\end{cor}

Note that Theorem \ref{thm1} does not directly apply to the torus
$(\T^{2n},\omega_L)$ with arbitrary linear symplectic forms
$\omega_L$ since there may not even exist such closed Lagrangian
submanifolds at all. The following result is proved using a similar
but somewhat different method.

\begin{theorem}\label{thm2}
The spectral metric on $\Ham(\T^{2n},\omega_L)$ does not extend to a
bi-invariant metric on $\Symp_0(\T^{2n},\omega_L)$ for the torus
$\T^{2n}$ with any linear symplectic form $\omega_L$.
\end{theorem}

One may conjecture as in the case for the Hofer norm that the above
theorems hold for general closed symplectic manifolds, although all
indications are that this problem is out of reach using currently
known methods.

\begin{guess}\label{conj1}
For any closed symplectic manifold $(M, \omega)$ such that
$Symp_0(M, \omega)$ is not identical to $\Ham(M, \omega)$, the
spectral metric on $\Ham(M, \omega)$ does not extend to a
bi-invariant metric on $\Symp_0(M, \omega)$.
\end{guess}

A related question is the bounded isometry conjecture proposed by
Lalonde and Polterovich in \cite{LP}. A symplectomorphism $\phi \in
\Symp_0(M,\omega)$ is said to be bounded with respect to the Hofer
norm $||\cdot||$ if
$$r_{||\cdot||}(\phi):=\sup_{f
\in \Ham(M,\omega)} \{||f \phi f^{-1} \phi^{-1}||\} < \infty.$$ Let
${\rm BI}_0(M,\omega, ||\cdot||)$ be the set of all bounded
symplectomorphisms with respect to the Hofer norm in
$\Symp_0(M,\omega)$. The bounded isometry conjecture for the Hofer
norm states that
$${\rm
BI}_0(M,\omega, ||\cdot||)=\Ham(M,\omega) \,\,\, \mbox{for  all}
\,\,\, (M,\omega).$$

As pointed out in \cite{Han2}, this is a very hard question in
general. Some partial results can be found in \cite{Han1, LP, LPe}.
Now the bounded isometry conjecture can easily be formulated with
respect to the spectral norm $\gamma$.

\begin{guess}\label{conj2}
${\rm BI}_0(M,\omega, \gamma)=\Ham(M,\omega)$ for all $(M,\omega)$.
\end{guess}

Both conjectures are trivial if $H^1(M, \R)$ vanishes. When $H^1(M,
\R)$ does not vanish, a necessary condition for both conjectures to
hold is that $\Ham(M, \omega)$ has infinite diameter with respect to
the spectral norm $\gamma$. This diameter is called the spectral
capacity of $M$ (cf. Equation (2.47) in Albers \cite{Alb}). As
noticed by Entov-Polterovich \cite{EnP} and explicitly pointed out
by McDuff \cite{Mc2}, the spectral capacity of certain manifolds
such as $\mathbb CP^n$ is finite. On the other hand, the spectral
capacity of $\T^{2n}$ is infinite. As far as I know, there are no
known examples of closed symplectic manifolds with finite spectral
capacity and such that $H^1(M, \R)$ does not vanish. Such examples
would immediately disprove both conjectures.

As one will see below, Conjecture \ref{conj2} implies Conjecture
\ref{conj1}. Note also that Conjecture \ref{conj2} implies the
corresponding conjecture for the Hofer norm, as the Hofer norm is
greater than the spectral norm. Therefore, it seems impossible to
prove them in general at this time. In this note, we only prove it
for the torus with linear symplectic forms.

\begin{theorem}\label{thm3}
The bounded isometry conjecture for the spectral norm holds for the
torus $\T^{2n}$ with all linear symplectic forms $\omega_L$.
\end{theorem}

\NI \textbf{Organization of the paper.} In Section \ref{tsn} we
review the definition of the spectral metric and its properties. We
prove Theorem \ref{thm1} in Section \ref{proofthm1} using an
explicit construction similar to that of Schwarz \cite{Sch}. Then we
prove Theorem \ref{thm3} in Section \ref{proofthm3}. Thereom
\ref{thm2} is a corollary of Theorem \ref{thm3}.

\bigskip

\NI \textbf{Acknowledgements.} The author wishes to thank Dusa
McDuff for her helpful comments on an earlier draft of this paper.

\section{The spectral metric}\label{tsn}
In this section we review without proofs the necessary ingredients
to define the spectral metric. The details can be found in Schwarz
\cite{Sch} for the case of closed symplectically aspherical
manifolds, and in Oh \cite{Oh1, Oh2} for closed symplectic manifolds
where $\omega$ is rational. See also Usher \cite{Ush} in the case
$\omega$ is not necessarily rational. For our purpose and for
simplicity, we only consider closed symplectically aspherical
manifolds, and will mainly follow Schwarz \cite{Sch} for the
definition. A symplectic manifold $(M, \omega)$ is symplectically
aspherical if
$${c_1}_{|\pi_2(M)}=\omega_{|\pi_2(M)}=0.$$

\subsection{Floer homology}
Let $(M, \omega)$ be a closed symplectically aspherical manifold.
Let $H_t=H_{t+1} : M \to \R$ be a time-dependent 1-periodic
Hamiltonian function. The time-dependent Hamiltonian vector field
${X_H}_t$ is given by $\iota({X_H}_t)=d H_t$. Consider the
Hamiltonian differential equation
\begin{equation} \label{eq1}
\dot{x}(t)={X_H}_t(x(t)).
\end{equation}
The solution of (\ref{eq1}) generates a family of Hamiltonian
diffeomorphisms $\phi_H^t : M \to M$ via
\begin{equation*}
\frac{d}{dt}{\phi}_H^t = {X_H}_t \circ \phi_t \,\,\, {\rm and}
\,\,\, \phi_0={\rm id}.
\end{equation*}
We shall first assume $H$ is regular, i.e. all fixed points of the
time-1 map $\phi_H^1$ are non-degenerate.

Given a time-dependent 1-periodic Hamiltonian function $H$, the
action functional $\A_H$ on the space $\mathcal{L}M$ of contractible
loops in $M$ is given by
$$\A_H(x) = \int_B u^*\omega + \int_0^1 H_t(x(t)) \, dt$$
for $x \in \mathcal{L}M$, where $B$ is the unit disk and $u : B \to
M$ is any extension of $x$. For symplectically aspherical manifold
$(M, \omega)$, the action functional $\A_H$ takes values in $\R$.
The critical points of $\A_H$ are contractible 1-periodic solutions
of (\ref{eq1}), i.e.
$${\rm Crit} \A_H = \mathcal{P}(H) : = \{x \in \mathcal{L}M \,|\, x
\,\,\, {\rm satisfies} \,\,\,(\ref{eq1}) \}.$$ The assumption that
$H$ is regular implies that $\mathcal{P}(H)$ is a finite set. The
Conley-Zehnder index $\mu_{\rm CZ} : \mathcal{P}(H) \to \mathbb{Z}$
takes values in $\mathbb{Z}$, and it can be normalized such that
\begin{equation*}
\mu_{\rm CZ}(x)=\mu_{\rm Morse}(x)
\end{equation*}
for $C^2$-small time-independent Morse function $H$, where $x \in
\mathcal{P}(H)= {\rm Crit} H$.

Let $J$ be an $\omega$-compatible almost complex structure on $M$.
We consider the moduli spaces $\M_{x, y}(J, H)$ of Floer
trajectories connecting $x, y \in \mathcal{P}(H)$. Namely,
\begin{align*}
\M_{x, y}(J, H):=\{u: \R \times S^1 \to M \, |& \, \partial_s u +
J(u) (\partial_t u - {X_H}_t(u))=0,\\& \lim_{s \to -\infty} u(s,
t)=x, \lim_{s \to + \infty} u(s, t)=y\}.
\end{align*}
Each solution $u \in \M_{x, y}(J, H)$ has finite energy
$$E(u):=\int_{-\infty}^{\infty} \int_0^1 |\partial_s u|^2 \, dsdt = \A_H(y)- A_H(x) \geq 0.$$
For generic $(J, H)$, $\M_{x, y}(J, H)$ is a smooth manifold of
dimension
$${\rm dim} \M_{x, y}(J, H)=\mu_{\rm CZ}(y)- \mu_{\rm CZ}(x).$$
When $\mu_{\rm CZ}(y)- \mu_{\rm CZ}(x)=1$, the space $\M_{x, y}(J,
H)$ is a 1-dimensional manifold, and the quotient $\M_{x, y}(J,
H)/\R$ is a set of finitely many isolated points. Here the
$\R$-action is the shift in $s$-variable. We denote $n(x, y):= \#
\M_{x, y}(J, H)/\R$, where the connecting orbits are counted with
appropriate signs.

Now consider the Floer chain complex $CF_*(H)=\mathcal{P}(H) \otimes
\Z$ with the grading given by the conley-Zehnder index $\mu_{\rm
CZ}$. The boundary operator is defined as
$$\partial : CF_*(H) \to CF_*(H),$$
$$\partial y = \sum_{\mu_{\rm CZ}(x)=\mu_{\rm CZ}(y)-1} n(x, y)\,
x.$$ By definition, the boundary operator $\partial$ has degree
$-1$.

Floer \cite{Fl} proved that $\partial \circ \partial =0$. Its
homology groups are called Floer homology groups of the pair $(J,
H)$, and are denoted by $HF_*(J, H)$. Floer also proved that
$HF_*(J, H)$ are independent of $(J, H)$, and are canonically
isomorphic to the singular homology groups of $M$.

\subsection{The PSS isomorphism} Let $f: M \to \R$ be a Morse
function, and $g$ be a generic Riemannian metric. Let $\rho : \R \to
[0,1]$ be any smooth cut-off function such that $\rho(s)=1$ for $s
\leq -1$, and $\rho(s)=0$ for $s \geq 1$. For $x \in \mathcal{P}(H)$
and $y \in {\rm Crit} f$, define the moduli spaces
\begin{align*}
\M^+_{x, y}(J, H; f, g):=\{& (u, \gamma) \,|\, u: \R \times S^1 \to
M, \,\, \gamma : [0, \infty) \to M,
\\& \partial_s u + J(u)
(\partial_t u - \rho(s){X_H}_t(u))=0, \dot{\gamma}+\nabla_g f \circ
\gamma=0,
\\& \int_{-\infty}^{\infty} \int_0^1 |\partial_s u|^2 \,
dsdt < \infty,
\\& \lim_{s \to -\infty} u(s, t)=x, \lim_{s \to
\infty} u(s, t)= \gamma(0), \lim_{s \to \infty} \gamma(s)=y\}.
\end{align*}
Similarly, define
\begin{align*}
\M^-_{y, x}(J, H; f, g):=\{&(\gamma, u) \, |\,  \gamma : (-\infty,
0] \to M, \,\, u: \R \times S^1 \to M,
\\& \dot{\gamma}+\nabla_g f \circ
\gamma=0, \partial_s u + J(u) (\partial_t u - \rho(-s){X_H}_t(u))=0,
\\& \int_{-\infty}^{\infty} \int_0^1 |\partial_s u|^2 \,
dsdt < \infty,
\\& \lim_{s \to -\infty} \gamma(s)=y, \gamma(0)=\lim_{s \to
-\infty} u(s, t), \lim_{s \to \infty} u(s, t)=x\}.
\end{align*}
These spaces are often called moduli spaces of $J$-holomorphic
spiked discs, where the spike is the gradient flow line between $y
\in {\rm Crit} f$ and $u(0)$. For generic choices of $(J, H, f, g)$,
the space $\M^+_{x, y}(J, H; f, g)$ is a manifold of dimension
$${\rm dim} \M^+_{x, y}(J, H; f, g) = \mu_{\rm Morse}(y)- \mu_{\rm
CZ}(x),$$ and the space $\M^-_{y, x}(J, H; f, g)$ is a manifold of
dimension
$${\rm dim} \M^-_{y, x}(J, H; f, g)= \mu_{\rm CZ}(x)- \mu_{\rm
Morse}(y).$$ Moreover, in dimension $0$ they are compact, hence
finite.

Define $\Phi_{\rm MF}$ from $CM_*(f)$ to $CF_*(H)$ by
$$\Phi_{\rm MF}(y)= \sum_{\mu_{\rm CZ}(x)= \mu_{\rm Morse}(y)} \#
\M^+_{x, y}(J, H; f, g) x,$$ and $\Phi_{FM}$ from $CF_*(H)$ to
$CM_*(f)$ by
$$\Phi_{\rm FM}(x)= \sum_{\mu_{\rm Morse}(y)= \mu_{\rm CZ}(x)} \#
\M^-_{y, x}(J, H; f, g) y.$$ Here the counting is again with
appropriate signs. These maps commute with boundary operators, hence
they induce homomorphisms
$$\Phi_{\rm MF} : HM_*(f, g) \to HF_*(J, H)$$
and
$$\Phi_{\rm FM} : HF_*(J, H) \to HM_*(f, g).$$
In fact they are isomorphisms, known as the PSS isomorphisms (cf.
\cite{PSS, Lu} for instance), and
$$\Phi_{\rm MF}= \Phi_{\rm FM}^{-1}.$$

\subsection{The spectral metric}
Given any $a \in \R$, define
$$\mathcal{P}^a(H)=\{x \in \mathcal{P}(H) \, | \, \A_H(x) \leq a\}$$
and
$$C_*^a(H)= \mathcal{P}^a(H) \otimes \Z.$$
Note that $C_*^a(H)$ is invariant under the Floer boundary operator
$\partial$. This means $(C_*^a(H),
\partial)$ is a sub-complex of $(C_*(H),
\partial)$, which is called filtered Floer chain complex. The
induced homology groups $HF_*^a(J, H)$ are called filtered Floer
homology groups.

The obvious inclusion map $$\iota_*^a : C_*^a(H) \to C_*(H)$$
induces $$\iota_*^a : HF_*^a(J, H) \to HF_*(J, H).$$ It is clear
that for $a$ sufficiently large, $HF_*^a(J, H) = HF_*(J, H)$ and
$\iota_*^a$ is the identity map; for $a$ sufficiently small,
$HF_*^a(J, H)=0$ and $\iota_*^a=0$.

Now given any nonzero homology class $\alpha \in H_*(M)$, under the
isomorphism $H_*(M) \cong HM_*(f, g)$, we can think of $\alpha$ as
an element in $HM_*(f, g)$. Define
$$C_{\alpha} (H) = \inf \,\{a \in \R \,|\, \Phi_{\rm MF} (\alpha) \in
{\rm im} \, \iota_*^a\}.$$ Note that $C_{\alpha} (H)$ are finite
numbers, and are critical values of the action functional $\A_H$.
Note also that they only depends on $H$, but not on $(f, g)$.

A standard energy estimate (cf. \cite{Sch}) shows that
\begin{equation} \label{eq2}
E_-(H) \leq C_{\alpha} (H) \leq E_+(H)
\end{equation}
where $$E_-(H)=\int_0^1 \, \min_{x \in M} H_t(x)\, dt, \,\,\,\,{\rm
and} \,\,\,\, E_+(H)=\int_0^1 \, \max_{x \in M} H_t(x) \, dt.$$ A
similar energy estimate shows that $C_{\alpha}$ satisfies
$$C_{\alpha} (H)- C_{\alpha} (K) \leq E_+(H-K).$$ This inequality
implies $C_\alpha$ is continuous with respect to the semi-norm
$$||H||= E_+(H)-E_-(H),$$ and can therefore be extended continuously
to non-regular Hamiltonians (Prop 2.14 in \cite{Sch}). Thus we do
not assume the regularity of $H$ any more from now on.

A deep result is that for a symplectically aspherical manifold $(M,
\omega)$,
$$C_{\alpha} (H)=C_{\alpha} (K) \,\,\,\, {\rm whenever} \,\,\,\,
\phi_H^1=\phi_K^1$$ where $\phi_H^1$ and $\phi_K^1$ are time-1 maps
of the Hamiltonians $H$ and $K$, respectively. Thus $C_\alpha$
descends to $\Ham(M, \omega)$ by defining $$C_\alpha(\phi)=
C_\alpha(H) \,\,\,\, {\rm if} \,\,\,\, \phi_H^1=\phi.$$ Note that in
general, $C_\alpha$ only descends to the universal cover
$\widetilde{\Ham}(M,\omega)$ (See \cite{Oh1, Oh2, Ush} for details).
In both cases, $C_\alpha$ are called the spectral invariants.

\begin{definition}
The spectral norm $\gamma : \Ham(M,\omega) \to \R$ is defined as
$$\gamma(\phi)= C_{[M]}(\phi) - C_{[pt]}(\phi)$$ where $[M] \in
H_{2n}(M)$ and $[pt] \in H_0(M)$ are the respective generators. The
induced metric
$$d_\gamma(\phi, \psi)= \gamma(\phi \psi^{-1})$$ is called the
spectral metric.
\end{definition}

\begin{remark} \rm
Since we are using homology $H_*(M)$ rather than cohomology
$H^*(M)$, our definition is slightly different from, but equivalent
to that of \cite{Sch}. We refer the reader to \cite{Oh1, Oh2, Ush}
for the definition of the spectral metric for general closed
symplectic manifolds.
\end{remark}

The following properties of the spectral norm $\gamma$ follow
directly from the properties of the spectral invariants $C_\alpha$.

\smallskip

(Nondegeneracy) $\gamma(\phi)\geq 0$, and $\gamma(\phi)=0
\Rightarrow \phi=id$.

(Symmetry) $\gamma(\phi)=\gamma(\phi^{-1})$.

(Triangle Inequality) $\gamma(\phi \psi)\leq
\gamma(\phi)+\gamma(\psi)$.

(Conjugate Invariancy) $\gamma(\psi \phi \psi^{-1})=\gamma(\phi)$.

\smallskip

These properties of $\gamma$ are equivalent to saying that the
spectral metric $d_\gamma$ is a nondegenerate bi-invariant metric.
In particular, $\gamma$ is conjugate-invariant implies $d_\gamma$ is
bi-invariant. Note that in view of (\ref{eq2}), the spectral norm
$\gamma$ is bounded from above by the Hofer norm $||\cdot||$. That
is,
$$\gamma(\phi) \leq ||\phi||.$$ Here the Hofer norm is another
conjugate-invariant norm on $\Ham(M, \omega)$ first discovered by
Hofer in \cite{Ho}. It is defined by
$$||\phi||=\inf_{\phi_H^1=\phi} \, \{\int_0^1 \, \max_{x \in M}
H_t(x)- \min_{x \in M} H_t(x) \, dt\}.$$ The induced metric
$$d_{\rm H}(\phi, \psi)= ||\phi \psi^{-1}||$$ is called the Hofer metric
which is also bi-invariant.

\section{Proof of Theorem \ref{thm1}}\label{proofthm1}
Recall from \cite{Han2} that to prove the spectral metric (in fact
any bi-invariant metric) $d_\gamma$ on $\Ham(M,\omega)$ does not
extend to a bi-invariant metric on $\Symp_0(M,\omega)$, it suffices
to show that there exists some unbounded symplectomorphism $\phi \in
\Symp_0(M, \omega)$ with respect to the spectral norm $\gamma$. That
is, $\phi$ satisfies
$$r_\gamma(\phi):=\sup_{f \in \Ham(M,\omega)} \{\gamma(f \phi
f^{-1} \phi^{-1})\} = \infty.$$ The reason is that assume $d_\gamma$
extend bi-invariantly, then we have
$$\gamma(f \phi f^{-1} \phi^{-1}) \leq \gamma(f \phi f^{-1}) +
\gamma(\phi^{-1})=\gamma(\phi)+\gamma(\phi^{-1}) = 2 \gamma(\phi).$$
This would imply $\gamma(\phi) \geq 2 r_\gamma(\phi) = \infty$,
which is a contradiction.

\bigskip

\NI \textbf{Proof of Theorem \ref{thm1}.} By assumption, there
exists some $\phi \in \Symp_0(M, \omega)$ such that $\phi(L) \cap L
= \emptyset$. It suffices to show that $\phi$ is unbounded with
respect to the spectral norm $\gamma$. To do this, one needs to find
a family of Hamiltonian diffeomorphisms $f_c$ such that $\gamma(f_c
\phi f_c^{-1} \phi^{-1}) \to \infty$ as $c \to \infty$. We follow
Schwarz (Example 5.7 in \cite{Sch}) for the construction of such
$f_c$'s.

The Lagrangian neighborhood theorem (cf. Weinstein \cite{Wei})
asserts that there exists a neighbourhood $V \subset M$ of $L$ which
is symplectomorphic to a neighbourhood $\mathcal{N}(L) \subset T^*L$
of the zero section. Now $L$ is a torus of dimension $n$, where $2n$
is the dimension of $M$, so $T^*L$ is symplectomorphic to $(\T^n
\times \R^n, \omega_0 = \sum_{i=1}^n\, dt_i \wedge dr_i)$, where
$t_i$'s and $r_i$'s are coordinates of $\T^n$ and $\R^n$
respectively. By shrinking $V$ if necessary, we can assume that
$\phi(V) \cap V = \emptyset$, and there exists a diffeomorphism
$\Phi : V \to \T^n \times B_\epsilon^n$ such that $\Phi^* \omega_0 =
\omega$. Here $B_\epsilon^n \subset \R^n$ is a ball of radius
$\epsilon$.

Given any constant $c$, let $\rho_c : \R^n \to \R$ be a smooth
function with $\mbox{supp}(\rho_c) \subset B_\epsilon^n$ and its
only critical values are $0$ and its maximum $\rho(\textbf{0})=c
>0$. Define $$\bar{\rho}_c: \T^n \times \R^n \to \R \,\,\, \mbox{such
that} \,\,\,\bar{\rho}_c(\textbf{t}, \textbf{r})
=\rho_c(\textbf{r}).$$ Using the diffeomorphism $\Phi : V \to \T^n
\times B_\epsilon^n$, we construct a family of time-independent
Hamitonians $F_c : M \to \R$ such that
$$F_c(x):=
\begin{cases}
\bar{\rho}_c \circ \Phi(x), & \mbox{if}\,\, x \in V,\\
\, 0, & \mbox{otherwise.}
\end{cases}$$
Let $f_c \in \Ham(M, \omega)$ be the time-$1$ map of the Hamiltonian
$F_c$. Since $f_c$ is supported in $V$ and $\phi(V) \cap V =
\emptyset$, it is easy to see that $$h_c := f_c \phi f_c^{-1}
\phi^{-1} \in \Ham(M, \omega)$$ is the time-$1$ map of the
Hamiltonian $H_c : M \to \R$ where
$$H_c(x):=
\begin{cases}
\bar{\rho}_c \circ \Phi(x), & \mbox{if}\,\, x \in V,\\
-\bar{\rho}_c \circ \Phi \circ \phi^{-1}(x), & \mbox{if}\,\, x \in \phi(V),\\
\,0, & \mbox{otherwise.}
\end{cases}$$

Note that since the Lagrangian torus $L$ is incompressible, the only
contractible $1$-periodic solutions of (\ref{eq1}) are the critical
points of $H_c$. Thus the action functional $\A_{H_c}$ has exactly
$3$ distinct critical values $c, 0, -c$. Recall that the spectral
norm $$\gamma(h_c)= C_{[M]} (h_c) - C_{[1]} (h_c)= C_{[M]} (H_c) -
C_{[1]} (H_c)$$ where $C_{[M]} (H_c)$ and $C_{[1]} (H_c)$ are two
distinct critical values of $\A_{H_c}$. So we get $\gamma(h_c) \geq
c$. Thus $\gamma(h_c) \to \infty$ as $c \to \infty$. This proves
$\phi$ is unbounded with respect to $\gamma$, hence Theorem
\ref{thm1}. \QED

\section{Proof of Theorem \ref{thm2} and \ref{thm3}}\label{proofthm3}

In this section, we prove the bounded isometry conjecture with
respect to the spectral norm for the torus $\T^{2n}$ with all linear
forms $\omega_L$, namely, Theorem \ref{thm3}. Theorem \ref{thm2}
then follows immediately as a corollary.

We begin with the flux homomorphism ${\rm Flux} : \widetilde{\rm
Symp}_0(M,\omega) \to H^1(M, \mathbb R)$ which is defined by $${\rm
Flux}(\{\phi_t\}):=\int_0^1 \iota(X_t)\omega \,dt$$ where the vector
filed $X_t$ is determined by $\frac{d}{dt} \phi_t=X_t \circ \phi_t$.
It induces
$${\rm Flux} : {\rm Symp}_0(M,\omega) \to H^1(M, \mathbb
R)/\Gamma_\omega$$ where $\Gamma_\omega:={\rm Flux}(\pi_1({\rm
Symp}_0(M,\omega)) \subset H^1(M, \mathbb R)$ is called the flux
subgroup, which is always discrete by Ono \cite{Ono}. It is well
known that the map ${\rm Flux}$ is surjective, and its kernel is
equal to $\Ham(M,\omega)$. Thus two symplectomorphisms have the same
flux if and only if they differ by a Hamiltonian diffeomorphism. See
\cite{MS1} Chapter 10 for more details.

Recall that $\phi \in \Symp_0(M,\omega)$ is bounded with respect to
$\gamma$ if $r_\gamma(\phi) < \infty$. Note that all Hamiltonian
diffeomorphisms $g$ are bounded since $r_\gamma(g) \leq 2 \gamma(g)<
\infty$. Similar to Proposition 1.2.A in \cite{LP}, $r_\gamma$
satisfies the triangle inequality $r_\gamma(\phi \psi) \leq
r_\gamma(\phi)+r_\gamma(\psi)$. Thus two symplectomorphisms with the
same flux are either both bounded or both unbounded.

From the above discussion, it is clear that to prove ${\rm
BI}_0(M,\omega, \gamma)=\Ham(M,\omega)$, it suffices to show that
for each nonzero value $v \in H^1(M, \R)/\Gamma_\omega$, there
exists some unbounded element $\phi \in \Symp_0(M,\omega)$ with
${\rm Flux}(\phi)=v$.

A linear symplectic form $\omega_L$ on $\T^{2n}$ is a $2$-form
$\omega_L=\sum_{i < j} \, a_{ij} dx_i \wedge dx_j$ such that
$\omega_L^n$ never vanishes. Let $\{\phi^i_\theta \} \in
\pi_1(\Symp_0(\T^{2n},\omega))$ be the loop of rotations of
$\T^{2n}$ along $x_i$ direction. One easily computes
$$\xi_i := {\rm Flux}(\{\phi^i_\theta\})=\displaystyle
\sum_{j=1}^{2n} a_{ij} dx_j.$$ Here we take the convention that
$a_{ij}=-a_{ji}$. In particular, $a_{ii}=0$. The following lemma can
be proved exactly the same way as Lemma 7.2 in \cite{Han1}.

\begin{lemma}
For $\T^{2n}$ with the linear symplectic form $\omega_L:=\sum_{i <
j}\, a_{ij} dx_i \wedge dx_j$, the flux subgroup $\Gamma_{\omega_L}
\subset H^1(\T^{2n}, \R)$ is generated by the above $\xi_i$'s over
$\Z$. That is, $\Gamma = \Z \langle \xi_1, \xi_2, \cdots, \xi_{2n}
\rangle$.
\end{lemma}

\begin{proof}
The proof is the same as the proof of Lemma 7.2 in \cite{Han1}, and
is therefore omitted.
\end{proof}

\smallskip

\NI {\bf Proof of Theorem \ref{thm3}.} Let $\phi \in
\Symp_0(\T^{2n},\omega_L)$ such that
$$\phi(x_1, x_2, \cdots, x_{2n})=(x_1+\alpha_1, x_2+\alpha_2, \cdots,
x_{2n}+\alpha_{2n})$$  where $\alpha_i \in \R/\Z$ for $1 \leq i \leq
2n$. Then
$${\rm Flux}(\phi)= \sum_{i=1}^{2n} \, \alpha_i \,
\xi_i.$$

As explained above, it suffices to show that $\phi$ is unbounded
with respect to the spectral norm $\gamma$ when at least one
$\alpha_i \in \R/\Z$ is nonzero. Assume $\alpha_1 \ne 0$ without
loss of generality. Thus $\phi (U) \cap U = \emptyset$ where $U
\subset \T^{2n}$ is defined by
$$U := \{(x_1, x_2, \cdots, x_4) \in \T^{2n} \mid |x_1|<
\epsilon \}.$$ for sufficiently small $\epsilon$.

As in the proof of Theorem \ref{thm1}, for every $c \in \R$, let
$F_c$ be a time-independent Hamiltonian function of $\T^{2n}$
supported in $U$ which has exactly two critical values, namely, $0$
and its maximum value $c$. We also require that $F_c$ depend only on
the first coordinate $x_1$. Denote by $f_c \in \Ham(\T^{2n},
\omega_L)$ the time-$1$ map of the Hamiltonian $F_c$. Since $\phi
(U) \cap U = \emptyset$, we know that $$h_c:= f_c \phi f_c^{-1}
\phi^{-1} \in \Ham(\T^{2n}, \omega_L)$$ is supported in the union of
two disjoint sets $U \cup \phi (U)$, and is the time-$1$ map of the
Hamiltonian $H_c=F_c-F_c \circ \phi^{-1}$, which has exactly three
critical values $c, 0, -c$. Using the fact that $\omega_L$ is a
linear symplectic form, one easily finds that the only contractible
$1$-periodic solutions of (\ref{eq1}) are the critical points of
$H_c$. Thus the spectral norm $\gamma(h_c)$ of $h_c$ goes to
infinity as $c$ goes to infinity, which implies $\phi$ is unbounded
with respect to $\gamma$. This proves Theorem \ref{thm3}. \QED

\smallskip

\begin{rmk} \rm
One may attempt to apply Theorem \ref{thm1} by showing that $\phi$
disjoins some Lagrangian torus $L \subset \T^{2n}$ from itself. For
a general linear symplectic form $\omega_L$, however, such
Lagrangian tori may not even exist in $\T^{2n}$. We go around the
difficulty by using the linearity of $\omega_L$.
\end{rmk}

\smallskip

\NI {\bf Proof of Theorem \ref{thm2}.} We have shown in the proof of
Theorem \ref{thm3} above the existence of some unbounded
symplectomorphism with respect to the spectral norm, which is
sufficient to prove Theorem \ref{thm2} as explained in Section
\ref{proofthm1}. \QED


\end{document}